\def\mystyle{}

\if\mystyle l
\documentclass[reqno,landscape]{amsart}
\else
\documentclass[reqno]{amsart}
\fi

\usepackage{begnac}

\newcommand{\brbinom}[2]{{#1 \brack #2}}

\begin{document}

\title{The Vandermonde determinant identity in higher dimension}

\author{Itaï \textsc{Ben Yaacov}}

\address{Itaï \textsc{Ben Yaacov} \\
  Université Claude Bernard -- Lyon 1 \\
  Institut Camille Jordan, CNRS UMR 5208 \\
  43 boulevard du 11 novembre 1918 \\
  69622 Villeurbanne Cedex \\
  France}

\urladdr{\url{http://math.univ-lyon1.fr/~begnac/}}

\thanks{Author supported by ANR project ValCoMo (ANR-13-BS01-0006) and ERC Grant no.\ 291111.}

% \date{\today}

\begin{abstract}
  We generalise the Vandermonde determinant identity to one which tests whether a family of hypersurfaces in $\bP^n$ has an unexpected intersection point.
\end{abstract}

\maketitle

\tableofcontents

\section*{Introduction}

The Vandermonde determinant identity tests by a single determinant whether a family of points on the line are distinct.
We generalise this to dimension $n$, testing by a single determinant whether some $n+1$ hyperplanes among a large family intersect.
This is further generalised for a family of hypersurfaces, up to an asymptotically negligible error.

\section{The linear identity}
\label{sec:LinearDualVandermonde}

We work in $\bP^n$.
Let $\fM_d$ denote the set of monomials of degree $d$ in $X = (X_0,\ldots,X_n)$.
We order it lexicographically, first comparing the exponents of $X_n$, if they are equal, comparing the exponent of $X_{n-1}$, and so on.
In particular, $\fM_1$ is ordered $X_0 < X_1 < \cdots < X_n$.

Given a finite set $I$, we let $\brbinom{I}{m} = \{J \subseteq I : |J| = m\}$, so $\left| \brbinom{I}{m} \right| = \binom{|I|}{m}$.
In addition, any $m \in \bN$ will be identified with the set $\{0, \ldots, m-1\}$, allowing us to write $\brbinom{m}{k}$.
If $I$ is an ordered set, we also equip $\brbinom{I}{m}$ with the lexicographic order: if $\sigma = \{s_0 < \cdots < s_{m-1}\}$ and $\tau = \{t_0 < \cdots < t_{m-1}\}$, then $\sigma < \tau$ if $s_{m-1} < t_{m-1}$, or if $s_{m-1} = t_{m-1}$ and $s_{m-2} < t_{m-2}$, and so on.
Passage to the complement reverses the order.

Let us also fix a number $m \in \bN$, which will pop up several times later.
It is a standard combinatorial fact that $|\fM_m| = \binom{m+n}{n} = \binom{m+n}{m}$.
If $(m+n) \setminus \sigma = \{s_0 < s_1 < \cdots < s_{n-1}\}$, then the unique order-preserving bijection between $\brbinom{m+n}{m}$ and $\fM_m$ sends $\sigma$ to
\begin{gather}
  \label{eq:MonomialFromSubset}
  \fm_\sigma
  = X_n^m \prod_{j<n} (X_j/X_{j+1})^{s_j-j}.
\end{gather}

Let $\Lambda = (\lambda_i : i \in I)$ be a family of linear forms with coefficients in a ring $A$, i.e., $\lambda_i \in A[X]_1$, and let $m = |I| - n$.
For each $\sigma \in \brbinom{I}{m}$, let
\begin{gather*}
  \Lambda^\sigma = \prod_{i \in \sigma} \lambda_i \in A[X]_m.
\end{gather*}
We define $M_\Lambda$ to be the square matrix of size $\binom{m+n}{n} \times \binom{m+n}{n}$, whose rows are the coefficient vectors of $\Lambda^\sigma$ for $\sigma \in \brbinom{I}{m}$, where columns correspond to monomials.
If we assume that $I$ is totally ordered, then both rows and columns of $M_\Lambda$ are ordered, so the sign of $\det M_\Lambda$ is well defined.

Let us consider some special cases, assuming for simplicity that $I = m+n$.
\begin{itemize}
\item If $m < 0$, i.e., if $|I| < n$, then everything still makes sense.
  In particular, the matrix $M_\Lambda$ is empty, so $\det M_\lambda = 1$.
\item If $m = 0$, then $M_\Lambda = (1)$ and again $\det M_\Lambda = 1$.
\item If $m = 1$, then $\brbinom{m+n}{m} = \brbinom{n+1}{1}$ is the set of singletons $\{i\}$ for $i \leq n$ with the natural order, $\fm_{\{i\}} = X_i$ and $\Lambda^{\{i\}} = \lambda_i$.
  In this case, $M_\Lambda$ is simply the matrix of coefficients of $\Lambda$, and we shall write $\det M_\Lambda = \det \Lambda$.
\end{itemize}

Now, to the general case.

\begin{prp}[Linear dual Vandermonde identity in arbitrary dimension]\label{prp:LinearDualVandermonde}
  With $\Lambda$ and $M_\Lambda$ as above, we have
  \begin{gather}
    \label{eq:LinearDualVandermonde}
    \det M_\Lambda = \prod_{\Lambda' \in \brbinom{\Lambda}{n+1}} \det \Lambda',
  \end{gather}
  where $\Lambda'$ is taken with the induced order from $\Lambda$.
\end{prp}
\begin{proof}
  This is an identity between polynomial expressions in the coefficients of $\Lambda$.
  We may therefore assume that all these coefficients are formal indeterminates, and that $A$ is the generated free $\bZ$-algebra.
  In particular, we may assume that $A$ is an integral domain, and then replace it with its fraction field $K$.

  A linear transformation of $K^{n+1}$, applied to $\Lambda$, can be decomposed into elementary transformations: substituting $X_j + a X_k$ ($k \neq j$) or $a X_j$ for $X_j$ (exchanging $X_j$ and $X_k$ can be obtained as a sequence of such).
  The first kind changes neither side of \autoref{eq:LinearDualVandermonde} (the effect on $M_\Lambda$ is that of a sequence of elementary transformations of the first kind applied directly to its rows).
  The second kind multiplies both by $a^N$, where $N = \binom{m+n}{n+1}$.
  Indeed, for the right hand side this is clear.
  For the left hand side, the product of all monomials of degree $m$ has total degree $m\binom{m+n}{n} = \frac{(m+n)!}{n!(m-1)!} = (n+1)N$, so the sum of degrees of $X_j$ in all monomials is $N$.
  Consequently, we may compose $\Lambda$ with any invertible linear transformation of $K^{n+1}$.

  We may assume that $I = m+n$, and proceed by induction on $(n,m)$.
  The case where $m \leq 0$ is covered by the special cases discussed earlier.
  If $n = 0$ then each $\lambda_i$ is of the form $a_i X_0$, and both sides of \autoref{eq:LinearDualVandermonde} are equal to $\prod a_i$.
  Therefore, we may assume that $m,n > 0$.
  By the preceding discussion, we may also assume that $\lambda_{m+n-1} = X_n$.

  Let us partition $\brbinom{m+n}{m}$ as $S_1 \cup S_2$, where
  \begin{gather*}
    S_1
    = \left\{ \sigma \in \brbinom{m+n}{m} : m+n-1 \notin \sigma \right\}
    = \brbinom{m+n-1}{m},
    \\
    S_2
    = \left\{ \sigma \in \brbinom{m+n}{m} : m+n-1 \in \sigma \right\}
    = \left\{ \sigma' \cup \{m+n-1\} : \sigma' \in \brbinom{m+n-1}{m-1} \right\}.
  \end{gather*}
  Since we are now changing $m$ and $n$, let us add them explicitly to the notation as superscripts.
  Thus, for example, if $\sigma \in S_1$, then $\fm^{m,n}_\sigma = \fm^{m,n-1}_\sigma$ (in which $X_n$ does not occur), and if $\sigma = \sigma' \cup \{m+n-1\} \in S_2$, then $\fm^{m,n}_\sigma = X_n \fm^{m-1,n}_{\sigma'}$.

  Let $\Lambda_1 = \bigl( \lambda_i(X_0,\ldots,X_{n-1},0) : i < m+n-1 \bigr)$, in dimension $n-1$, with indeterminates $(X_0,\ldots,X_{n-1})$, and $\Lambda_2 = \bigl( \lambda_i : i < m+n-1 \bigr)$.
  If $\sigma \in S_1$, then $\Lambda^\sigma = \Lambda_1^\sigma + X_n P$ for some polynomial $P$, and if $\sigma = \sigma' \cup \{m+n-1\} \in S_2$, then $\Lambda^\sigma = \lambda_{m+n-1} \Lambda^{\sigma'} = X_n \Lambda_2^{\sigma'}$.
  Therefore,
  \begin{gather*}
    M^{m,n}_\Lambda =
    \begin{pmatrix}
      M^{m,n-1}_{\Lambda_1} & ?
      \\
      0 & M^{m-1,n}_{\Lambda_2}
    \end{pmatrix}.
  \end{gather*}
  By our induction hypotheses,
  \begin{gather*}
    \det M^{m,n-1}_{\Lambda_1}
    = \prod_{\Lambda' \in \brbinom{\Lambda_1}{n}} \det \Lambda'
    = \prod_{\Lambda' \in \brbinom{\Lambda_2}{n}} \det (\Lambda', X_n)
    = \prod_{\Lambda' \in \brbinom{\Lambda_2}{n}} \det (\Lambda', \lambda_{m+n-1}),
    \\
    \det M^{m-1,n}_{\Lambda_2}
    = \prod_{\Lambda' \in \brbinom{\Lambda_2}{n+1}} \det \Lambda'.
  \end{gather*}
  Our assertion follows.
\end{proof}

Let us consider the case where $n = 1$, with the indeterminates $X,Y$.
Let $\Lambda = (\lambda_i : i \leq m)$, say $\lambda_i = - b_i X + a_i Y$.
We have $\brbinom{m+1}{m} = \bigl\{ (m+1) \setminus \{i\} : i \leq m \bigr\}$, which we may identify with the set $m+1$ with reversed order.
Since this affects both rows and columns, we may replace it with the usual order.

For $i \leq m$, the monomial $\fm_{(m+1) \setminus \{i\}}$ is $X^i Y^{m-i}$.
Accordingly, let $M'_\Lambda$ be the matrix $(a_j^i b_j^{m-i})_{i,j \leq m}$: its $j$th column is the image of $(a_j,b_j)$ under the $m$-fold Veronese map.
The $(i,j)$ entry of $M_\Lambda M'_\Lambda$ is $\Lambda^{(m+1) \setminus \{i\}}(a_j,b_j) = \prod_{k \neq i} (a_k b_j - a_j b_k)$.
This vanishes if $i \neq j$, so
\begin{gather*}
  \det M_\Lambda \det M'_\Lambda
  = \prod_{j \neq i} \, (a_j b_i - a_i b_j).
\end{gather*}
By \autoref{prp:LinearDualVandermonde}, $\det M_\Lambda = \prod_{i < j} \, (a_i b_j - a_j b_i)$, so
\begin{gather}
  \label{eq:LinearVandermonde1}
  \det M'_\Lambda
  = \prod_{i < j} \, (a_j b_i - a_i b_j).
\end{gather}
Substituting $b_i = 1$, we get the classical (affine, dimension one) Vandermonde identity, which in turn easily implies \autoref{eq:LinearVandermonde1}:
\begin{gather}
  \label{eq:AffineVandermonde}
  \det
  \begin{pmatrix}
    1 & a_0 & \cdots & a_0^m \\
    1 & a_1 & \cdots & a_1^m \\
    \vdots & \vdots & & \vdots \\
    1 & a_m & \cdots & a_m^m \\
  \end{pmatrix}
  = \prod_{i<j\leq m} (a_j - a_i).
\end{gather}

In other words, for $n=1$, the identities \autoref{eq:AffineVandermonde} and \autoref{eq:LinearDualVandermonde} are equivalent, and in a sense dual, one referring to linear forms and the other to their zeros.

For $n > 1$, one can still generalise the classical Vandermonde identity as a dual of \autoref{eq:LinearDualVandermonde}, but the construction of $M'_\Lambda$ is more tedious.
For any $\sigma = (m+n) \setminus \{s_0 < \cdots < s_{n-1}\} \in \brbinom{m+n}{m}$, let $a_\sigma$ be the wedge product $\lambda_{s_0} \wedge \ldots \wedge \lambda_{s_{n-1}}$, namely the vector of minors (with appropriate signs) of the $(n+1) \times n$ matrix $(\lambda_{s_0} \cdots \lambda_{s_{n-1}})$.
Let $V_m(a_\sigma)$ be its image under the $m$-fold Veronese map, namely the vector of values, at $a_\sigma$, of all monomials of degree $m$.
We define $M'_\Lambda$ to be the $\brbinom{m+n}{m} \times \brbinom{m+n}{m}$ matrix whose $\sigma$th column is the vector $V_m(a_\sigma)$, with our usual indexing of monomials by $\brbinom{m+n}{m}$.

\begin{cor}[Linear Vandermonde identity in arbitrary dimension]
  \label{cor:LinearVandermonde}
  With $\Lambda$ and $M'_\Lambda$ as above, we have
  \begin{gather}
    \label{eq:LinearVandermonde}
    \det M'_\Lambda
    = \prod_{\Lambda' \in \brbinom{\Lambda}{n+1}} (-1)^{\binom{n+1}{2}} (\det \Lambda')^n,
  \end{gather}
  where $\Lambda'$ is taken with the induced order from $\Lambda$.
\end{cor}
\begin{proof}
  The argument follows the same path as for $n = 1$.
  The $(\sigma,\tau)$ entry of $M_\Lambda M'_\Lambda$ is $\Lambda^\sigma(a_\tau) = \prod_{k \in \sigma} \det \bigl( \lambda_k, (\lambda_j)_{j \notin \tau} \bigr)$, so again the product matrix is diagonal, and
  \begin{gather}
    \label{eq:LinearVandermondeDuality}
    \det M_\Lambda \det M'_\Lambda = \prod_{\Lambda' \in \brbinom{\Lambda}{n+1}} (-1)^{\binom{n+1}{2}} (\det \Lambda')^{n+1}.
  \end{gather}
  Indeed, each $\det \Lambda'$ occurs $n+1$ times as $\det \bigl( \lambda_k, (\lambda_j)_{j \notin \tau} \bigr)$, once for each choice of $\lambda_k \in \Lambda'$, giving a total of $0 + 1 + \cdots + n = \binom{n+1}{2}$ changes of sign.
  Now divide by \autoref{eq:LinearDualVandermonde}.
\end{proof}

We observe that in dimension $n > 1$, the statement of the dual identity is more elegant.
It is also the one that admits a relatively simple ``asymptotic'' non-linear generalisation, in which that $\lambda_i$ are replaced with homogeneous polynomials of some higher degree.

\begin{rmk}[Ulysse Serres's « usine à points »]
  It follows from \autoref{cor:LinearVandermonde} that given a family of $n+d$ hyperplanes in general position in $n$-dimensional space, no polynomial of degree $d$ can vanish at all intersection points of subfamilies of $n$ many hyperplanes.
  The classical Vandermonde determinant identity \autoref{eq:AffineVandermonde} provides us with a particularly easy way to produce such a family of hyperplanes and calculate their intersection points and an analogue of \autoref{eq:LinearVandermondeDuality} allows us to prove that the determinant of the matrix of monomials is non-zero.
  We obtain the following self-contained argument.

  Let $n,d \in \bN$.
  Let $I = \{1,\ldots,n+d\}$, and assume that $a_1,\ldots,a_{n+d}$ are distinct members of some field $k$.

  Let $\Sigma$ denote the collection of all subsets of $I$ of size $n$.
  For $\sigma \in \Sigma$, define a polynomial
  \begin{gather*}
    f_\sigma(T) = \prod_{i \in \sigma} \, (T - a_i) = T^n + b_{\sigma,1} T^{n-1} + \cdots + b_{\sigma,n}.
  \end{gather*}
  Let $b_\sigma = (b_{\sigma,1},\ldots,b_{\sigma,n}) \in k^n$ be the vector of coefficients.

  Let $X = (X_1,\ldots,X_n)$ denote an indeterminate point in $k^n$.
  For $1 \leq i \leq n+d$, define an affine function $\lambda_i\colon k^n \rightarrow k$ by
  \begin{gather*}
    \lambda_i(X) = a_i^n + a_i^{n-1} X_1 + \cdots + X_n.
  \end{gather*}
  Then
  \begin{gather*}
    \lambda_i(b_\sigma) = a_i^n + a_i^{n-1} b_{\sigma,1} + \cdots + b_{\sigma,n} = f_\sigma(a_i).
  \end{gather*}

  Let $\Xi = \{\xi \in \bN^n : \xi_1 + \cdots + \xi_n \leq d\}$.
  If $\sigma \in \Sigma$, then $|I \setminus \sigma| = d$, and we may define a polynomial
  \begin{gather*}
    g_\sigma(X) = \prod_{i \in I \setminus \sigma} \lambda_i(X) = \sum_{\xi \in \Xi} c_{\sigma,\xi} X^\xi.
  \end{gather*}

  Let us define two matrices $B = (b_\sigma^\xi)$ and $C = (c_{\sigma,\xi})$, whose rows are indexed by $\xi \in \Xi$ and columns by $\sigma \in \Sigma$.
  Let $D = C^t \, B = (d_{\sigma,\tau})$, where $\sigma,\tau \in \Sigma$.
  Then
  \begin{gather*}
    d_{\sigma,\tau}
    = \sum_{\xi \in \Xi} c_{\sigma,\xi} b_\tau^\xi
    = g_\sigma(b_\tau)
    = \prod_{i \in I \setminus \sigma} \lambda_i(b_\tau)
    = \prod_{i \in I \setminus \sigma} f_\tau(a_i)
    = \prod_{i \in I \setminus \sigma} \prod_{j \in \tau} \, (a_j - a_i).
  \end{gather*}
  Therefore, $d_{\sigma,\tau} \neq 0 \ \Longleftrightarrow  \ \tau \cap (I \setminus \sigma) = \emptyset \ \Longleftrightarrow \ \tau \subseteq \sigma$.
  Since $|\tau| = |\sigma| = n$, we have $d_{\sigma,\tau} \neq 0$ if and only if $\sigma = \tau$.
  Therefore $D$ is diagonal, and $\det D \neq 0$.

  It is a standard combinatorial fact that $|\Xi| = \binom{n+d}{n}$, so $B$ and $C$ are square matrices, and since $\det D \neq 0$, they are invertible.
  Let $h = \sum_{\xi \in \Xi} h_\xi X^\xi$ be a non-zero polynomial of degree at most $d$, and let $H = (h_\xi : \xi \in \Xi)$.
  Since $B$ is invertible, $H B = \bigl( h(b_\sigma) : \sigma \in \Sigma \bigr)$ is non-zero.
  We conclude that the only polynomial of degree at most $d$ that vanishes at all the points $b_\sigma$ is the zero polynomial.

  Finally, by \autoref{eq:LinearDualVandermonde} and \autoref{eq:AffineVandermonde}, we have
  \begin{gather*}
    \det C = \prod_{i<j} \prod_{\{i,j\} \subseteq \sigma \in \brbinom{n+d}{n+1}} \, (a_j-a_i)
    = \prod_{i<j} \, (a_j-a_i)^{\binom{n+d-2}{n-1}}.
  \end{gather*}
  By our calculations,
  \begin{gather*}
    \det D = \prod_{\sigma \in \Sigma} \prod_{i \in I \setminus \sigma} \prod_{j \in \sigma} \, (a_j-a_i)
    = \pm \prod_{i<j} \, (a_j-a_i)^{2 \binom{n+d-2}{n-1}} = \pm \det C^2.
  \end{gather*}
  Therefore $\det B = \pm \det C$.
\end{rmk}

\section{An asymptotic higher-degree identity}
\label{sec:NonLinearDualVandermonde}

In \autoref{prp:LinearDualVandermonde}, a line of the matrix $M_\Lambda$ represent a polynomial $\Lambda^\sigma$, which may be viewed as a split polynomial $f = \prod_{i<m+n} \lambda_i$ from which $n$ linear factors are dropped.
In dimension $n > 1$, polynomials do not necessarily split.
This introduces two related issues.
First of all, we need to replace the determinant of a family of $n+1$ linear forms, with the more general notion of resultant.

\begin{fct}
  \label{fct:MacaulayResultant}
  Let $n \in \bN$, and let $X$ denote a tuple of $n+1$ indeterminates, representing homogeneous coordinates in $\bP^n$.
  For each $i \leq n$, let $d_i \geq 1$ and let $T^*_i$ be a tuple of indeterminates representing the coefficients of a polynomial in $X$, homogeneous of degree $d_i$.
  Then there exists an irreducible polynomial
  \begin{gather*}
    \fR \in \bZ[T^*_i : i \leq n]
  \end{gather*}
  called the \emph{Macaulay resultant}, with the following property:
  \begin{quote}
    For any algebraically closed field $K$, and any family of homogeneous polynomials $f_i \in K[X]_{d_i}$, for $i \leq n$, we have $\fR(f_{\leq n}) = 0$ if and only if the family $f_{\leq n}$ has a common zero in $\bP^n(K)$.
  \end{quote}
  This polynomial is unique up to sign.
  It is homogeneous in each $T^*_i$, of degree $\prod_{j \neq i} d_j$.
\end{fct}
\begin{proof}
  Macaulay \cite{Macaulay:Elimination}.
\end{proof}

We may even fix the sign of $\fR$, as follows.
When $d_0 = d_1 = \cdots = d_n = 1$, we have $\fR(\lambda_0,\ldots,\lambda_n) = \pm \det(\lambda_0,\ldots,\lambda_n)$, where $\lambda_i$ are linear forms and we take the determinant of the coefficient matrix.
In particular, $\fR(X_0,\ldots,X_n) = \pm 1$, and we may choose the sign so $\fR(X_0,\ldots,X_n) = 1$.
In higher degrees, we have $\fR(\lambda_0^{d_0},\ldots,\lambda_n^{d_n}) = \pm \det(\lambda_0,\ldots,\lambda_n)^{\prod d_i}$, and again, we may choose the sign so $\fR(X_0^{d_0},\ldots,X_n^{d_n}) = 1$.
Once the signs have been fixed in this fashion, we have $\fR(\ldots, fg, \ldots) = \fR(\ldots, f, \ldots) \fR(\ldots, g, \ldots)$.

We may extend the definition in a coherent fashion to the case where $d_i = 0$ for some $i$, letting $\fR(\ldots,T^*_i,\ldots) = (T^*_i)^{\prod_{j \neq i} d_j}$ (though this is not irreducible).
If two degrees vanish, then $\fR = 1$.

When $n = 1$, $\fR$ takes the familiar form as the determinant of the Sylvester matrix, with entries taken from the coefficients $T^*_0$ and $T^*_1$.

\begin{ntn}
  \label{ntn:MacaulayResultant}
  Since the resultant corresponds to intersections, we shall express $\fR(f_0,\ldots,f_n)$ as $f_0 \wedge \cdots \wedge f_n$.
  When $F = (f_0,\ldots,f_n)$, we are going to write $\fR(F)$ as $F^\wedge$.
\end{ntn}

The second issue is that and dropping a single linear factor from a polynomial that does not split no longer makes sense.

\begin{dfn}
  \label{dfn:PseudoSplitPolynomial}
  A \emph{pseudo-split polynomial} $f$ over $A$, in indeterminates $X$ and of degree $m$, consists of a set $I_f$ of (indices of) \emph{pseudo-factors}, together with polynomials $f^\sigma \in A[X]_{|\sigma|}$ for each $\sigma \subseteq I_f$, where $f^\emptyset = 1$.
  \begin{itemize}
  \item If $\sigma \subseteq I_f$, then $g = f^\sigma$ may be viewed as a pseudo-split polynomial, with $I_g = \sigma$ and $g^\tau = f^\tau$ for $\tau \subseteq \sigma$.
  \item If $f,g$ are pseudo-split and $I_f \cap I_g = \emptyset$, then $fg$ is pseudo-split with $I_{fg} = I_f \cup I_g$ and $(fg)^\sigma = f^{\sigma \cap I_f} g^{\sigma \cap I_g}$.
  \end{itemize}
\end{dfn}

Intuitively, we think of $f^\sigma$ as the product of the pseudo-factors indexed by $\sigma$.
Notice, however, that we do \emph{not} require that $f^{\sigma \cup \tau} = f^\sigma f^\tau$.
When we want to view a pseudo-split polynomial $f$ as a plain polynomial, we identify it with $f^{I_f}$ (all pseudo-factors retained).

\begin{dfn}
  \label{dfn:PseudoSplitMF}
  Let $f$ be pseudo-split of degree $m+n$, and assume that $I = I_f$ is ordered.
  We define $M_f$ to be the matrix whose rows are the coefficient vector of $f^\sigma \in A[X]_m$ for $\sigma \in \brbinom{I}{m}$.
  We order $\brbinom{I}{m}$ and the monomials $\fM_m$ as in \autoref{sec:LinearDualVandermonde}.

  Given a family $F = (f_i : i < k)$ of pseudo-split polynomials with disjoint ordered $I_{f_i}$, we define $I_F = I_{\prod F} = \bigcup I_{f_i}$, ordered so $I_{f_i} < I_{f_{i+1}}$, and $M_F = M_{\prod f_i}$.
\end{dfn}

A split polynomial $f = \prod_{i \in I} \lambda_i$ is pseudo-split in the obvious fashion, namely, $f^\sigma = \prod_{i \in \sigma} \lambda_i$.
Let $\Lambda = (\lambda_i : i < m+n)$ be a family of linear forms and $f = \prod \Lambda$ a pseudo-split polynomial in this sense.
Then $I_f = m+n$, and $f^\sigma$ is exactly $\Lambda^\sigma$ of \autoref{sec:LinearDualVandermonde}.
It follows that $M_\Lambda$, as per \autoref{dfn:PseudoSplitMF}, coincides with the one defined in \autoref{sec:LinearDualVandermonde}, which we are attempting to generalise.

At the other extreme, an \emph{indeterminate} pseudo-split polynomial of degree $m+n$ is one where the coefficients of $f^\sigma$ (with the exception of $f^\emptyset = 1$) are indeterminates over the base ring.
For such a polynomial, $M_f$ is an indeterminate matrix (as soon as $m \geq 1$), and we cannot expect any useful identities involving its determinant.

We are interested in a situation that lies somewhere in between these two extremes.
We are given a polynomial $f$ that factors as $\prod_{i<k} f_i$, where the family $F = (f_i : i < k)$ consists of polynomials of relatively low degrees $d_i = \deg f_i$, and $m+n = \deg f = \sum d_i$.
With no additional information being provided, we might as well take the $f_i$ to be indeterminate pseudo-split polynomials, with $I_{f_i}$ ordered and pairwise disjoint.
We propose to prove an expression for $\det M_F$, generalising \autoref{prp:LinearDualVandermonde}, that is asymptotically useful as $k$ grows, provided that the degrees $d_i$ remain bounded.
Let $T^*$ denote the indeterminate coefficients, so all the polynomials involved are in $A[X]$ where $A = \bZ[T^*]$.

\begin{lem}
  \label{lem:VandermondeHypersurfacesSpecialisation}
  Let $\Lambda = (\lambda_j : j \in I_F)$ be indeterminate linear forms, with coefficients $S^*$, and define a morphism
  \begin{gather*}
    \varphi\colon \bZ[T^*] \rightarrow \bZ[S^*]
  \end{gather*}
  by sending each $f_i$ to $\prod_{j \in I_{f_i}} \lambda_j$ as a pseudo-split polynomial, namely, sending each coefficient of $f_i^\sigma$ (which is one of the generators $T^*$) to the corresponding coefficient of $\varphi(f_i^\sigma) = \Lambda^\sigma$.
  \begin{enumerate}
  \item
    \label{item:VandermondeHypersurfacesSpecialisationMF}
    With $M_F$ defined as per \autoref{dfn:PseudoSplitMF},
    \begin{gather*}
      \varphi(\det M_F) = \prod_{\Lambda' \in \brbinom{\Lambda}{n+1}} \det \Lambda'.
    \end{gather*}
    Consequently, $\det M_F \neq 0$, and has no non-zero factors in $\bZ$.
  \item
    \label{item:VandermondeHypersurfacesSpecialisationEquality}
    If $P,Q \in \bZ[T^*]$ are homogeneous, $\varphi(P) \neq 0$, $P \mid Q$ and $\varphi(Q) \mid \varphi(P)$, then $P = \pm Q$.
    If $\varphi(P) = \varphi(Q)$, then $P = Q$.
  \end{enumerate}
\end{lem}
\begin{proof}
  The first assertion is just \autoref{prp:LinearDualVandermonde}.
  For the second, we have $PR = Q$, with $R$ homogeneous, and $\varphi(P) \varphi(R) \mid \varphi(P)$.
  Therefore $\varphi(R) = \pm 1$, which is only possible if $R = \pm 1$.
  If $\varphi(P) = \varphi(Q)$, then $R = 1$.
\end{proof}

\begin{lem}
  \label{lem:MGCD}
  Let $G,H \subseteq F$.
  Then $\det M_G \mid \det M_F$, and (up to sign)
  \begin{gather*}
    \det M_{G \cap H} = \gcd(\det M_G,\det M_H).
  \end{gather*}
\end{lem}
\begin{proof}
  Let $h = \prod (F \setminus G)$.
  If $\sigma \subseteq I_G$, then $F^{\sigma \cup I_h} = G^\sigma h$.
  Assuming that $I_G$ is an final segment, sets of the form $\sigma \cup I_h$ form an initial segment, and we have
  \begin{gather*}
    M_F =
    \begin{pmatrix}
      M_G A \\ B
    \end{pmatrix}
    =
    \begin{pmatrix}
      M_G  & 0 \\ 0 & I
    \end{pmatrix}
    \begin{pmatrix}
      A \\ B
    \end{pmatrix},
  \end{gather*}
  where $A$ is the matrix of multiplication by $h$.
  Therefore $\det M_G \mid \det M_F$.

  In particular, if $G,H \subseteq F$, then $\det M_{G \cap H} \mid \gcd(\det M_G,\det M_H)$.
  With $\varphi$ as in \autoref{lem:VandermondeHypersurfacesSpecialisation}, we have $\varphi(\det M_{G \cap H}) = \gcd\bigl( \varphi(\det M_G), \varphi(\det M_H) \bigr)$.
  The desired equality (up to sign) now follows from \autoref{lem:VandermondeHypersurfacesSpecialisation}\autoref{item:VandermondeHypersurfacesSpecialisationEquality}.
\end{proof}

Let us define
\begin{gather*}
  \Delta_F = \prod_{G \subseteq F} \det M_G^{(-1)^{|F \setminus G|}},
  \qquad
  \text{so}
  \qquad
  \det M_F = \prod_{G \subseteq F} \Delta_G.
\end{gather*}
This inclusion/exclusion stratification applies of course for any map from finite sets into an Abelian (in our case, multiplicative) group.
A more sophisticated inclusion/exclusion principle allows us to calculate the least common multiple in a unique factorisation domain (or more generally, the supremum in a lattice-ordered Abelian group):
\begin{gather*}
  \lcm(P_j : j \in J) = \prod_{\emptyset \neq K \subseteq J} \gcd(P_k : k \in K)^{(-1)^{|K|-1}}.
\end{gather*}
Let us apply this to $P_j = \det M_{G_j}$, for $G_j \subseteq F$.
Let $G_K = \bigcap_{k \in K} G_k$.
Then by \autoref{lem:MGCD},
\begin{gather*}
  \lcm(\det M_{G_j} : j \in J) = \prod_{\emptyset \neq K \subseteq J} \det M_{G_K}^{(-1)^{|K|-1}}.
\end{gather*}
Applying this to the family of $F \setminus\{f\}$, for $f \in F$, we obtain
\begin{gather*}
  \lcm\left( \det M_{F \setminus \{f\}} : f \in F \right)
  = \prod_{\emptyset \neq K \subseteq F} \det M_{F \setminus K}^{(-1)^{|K|-1}}
  = \det M_F / \Delta_F.
\end{gather*}
Since $\lcm\left( \det M_{F \setminus \{f\}} : f \in F \right)$ divides $\det M_F$, the expression $\Delta_F$ is a polynomial in the coefficients of $F$.
In particular, it can be evaluated for any family of pseudo-split polynomials, over any ring.

Let us attempt to calculate $\Delta_F$ explicitly, in three cases, which we dub ``easy'', ``intractable'', and ``interesting''.
First, the easy case, when $|F| > n + 1$.
Applying \autoref{lem:VandermondeHypersurfacesSpecialisation}, we have $\varphi(\Delta_F) = 1$ (and $1 \mid \Delta_F$), so $\Delta_F = 1$.
Next, the intractable case, when $|F| \leq n$.
In this case we cannot improve our prior observation that $\Delta_F$ is a polynomial (and not a rational function, as its form may suggest at first).

Finally, the interesting case, when $|F| = n+1$.
If, under some specialisation of $\bZ[T^*]$ into an algebraically closed field $K$, the family $F$ has a common zero in $\bP^n(K)$, then this is a zero of $F^\sigma$ for all $\sigma \in \brbinom{I_F}{m}$ (since $\sigma$ must contain at least one $I_{f_i}$ entirely), and the matrix $M_F$ specialises to a singular one.
In other words, every zero of $F^\wedge \in \bZ[T^*]$ (the resultant, see \autoref{ntn:MacaulayResultant}) is a zero of $\det M_F$.
Since $F^\wedge$ is irreducible, it divides $\det M_F$, and cannot divide $\det M_G$ for any $G \subsetneq F$, so $F^\wedge \mid \Delta_F$.
Let us evoke \autoref{lem:VandermondeHypersurfacesSpecialisation} once more.
Since the resultant is multiplicative, and since the resultant of $n+1$ linear forms is their determinant, we have
\begin{gather*}
  \varphi(F^\wedge) = \prod \bigl\{ \det (\lambda_{i_0}, \ldots, \lambda_{i_n}) : i_0 \in I_{f_0}, \ldots, i_n \in I_{f_n}\bigr\}.
\end{gather*}
It is easy to check that $\varphi(\Delta_F)$ is equal to the same, so $\Delta_F = F^\wedge$.
(This is valid even when $d_i = 0$ for some $i$, in which case $f_i = 1$ and $\Delta_F = F^\wedge = 1$).

We have thus proved:

\begin{thm}[The dual Vandermonde identity in arbitrary dimension and degree]
  \label{thm:NonLinearDualVandermonde}
  For any family $F = (f_i : i < k)$ of pseudo-split polynomials:
  \begin{gather}
    \label{eq:NonLinearDualVandermonde}
    \det M_F = \prod_{H \in \brbinom{F}{n+1}} H^\wedge \prod_{G \subseteq F, \ |G| \leq n} \Delta_G.
  \end{gather}
\end{thm}

Assume that we are given a family $F$ of (ordinary) polynomials.
In order to apply \autoref{thm:NonLinearDualVandermonde}, we need to make them pseudo-split, so we introduce, in an arbitrary fashion, auxiliary polynomials $f_i^\sigma$ for $\emptyset \neq \sigma \subsetneq I_{f_i}$.
The product of resultants on the right hand side of \autoref{eq:NonLinearDualVandermonde} is the ``main term'', which only depends on the family $F$, and is analogous to the right hand side of \autoref{eq:LinearDualVandermonde}.
The product of $\Delta_G$ is the ``error term'', which may also depend on the auxiliary coefficients, and has no obvious meaning.

In the special case where $\deg f_i = 1$ for all $i$, there are no auxiliary coefficients nor error term, and \autoref{thm:NonLinearDualVandermonde} specialises to \autoref{prp:LinearDualVandermonde} (the $\Delta_G$ are all equal to one, since $M_G$ is $I_1$ for $|G| = n$, and $I_0$ for $|G| < n$).
In the general case, the error term is present.
However, as $F$ grows, the number of $\Delta_G$ grows like $|F|^n$, while the number of $H^\wedge$ grows like $|F|^{n+1}$.
Therefore, the error term is negligible \emph{asymptotically}.

From \autoref{thm:NonLinearDualVandermonde} we may recover a \emph{relative} identity with respect to a hypersurface.
Again, this is an asymptotic result, with growth rates $|F|^{n-1}$ for the error term, against $|F|^n$ for the main term.

\begin{cor}
  \label{cor:NonLinearDualVandermondeRelative}
  Let $F,g$ denote the sequence $F$ followed by $g$, and assume that $M_F$ is nonsingular.
  Then
  \begin{gather}
    \label{eq:NonLinearDualVandermondeRelative}
    \frac{\det M_{F,g}}{\det M_F} = \prod_{H \in \brbinom{F}{n}} (H,g)^\wedge \prod_{G \subseteq F, \ |G| < n} \Delta_{G,g}.
  \end{gather}
\end{cor}

The usefulness of either identity is of course conditioned on the non-vanishing of any of the terms.
Let us conclude with some particularly easy sufficient conditions for that.

\begin{lem}
  \label{lem:NonLinearDualVandermondeDeltaNonZero}
  Let $A$ be an integral domain, and let $F$ be a finite family of pseudo-split polynomials, generic over $A$.
  In other words, $F$ is a family of pseudo-split polynomials over some larger ring $B$, that is freely generated over $A$ by the coefficients of the polynomials $f^\sigma$, where $f \in F$ and $\emptyset \neq \sigma \subseteq I_f$.
  Then $\det M_F \neq 0$ and $\Delta_F \neq 0$ (when evaluated in $B$).

  Moreover, assume that $g \neq 0$ is an additional pseudo-split polynomial, such that $(g^\sigma : \emptyset \subsetneq \sigma  \subsetneq I_g)$ are also generic (but not necessarily $g$ itself).
  Then $\det M_{F,g} \neq 0$ and $\Delta_{F,g} \neq 0$.
\end{lem}
\begin{proof}
  It is enough to show that $M_F \neq 0$ and $M_{F,g} \neq 0$, and in order to show this, we may specialise the members of $F$ to split polynomials with generic linear factors.
  In other words, we may assume that $F$ is a family of generic linear forms.
  Then $\det M_F \neq 0$ by \autoref{prp:LinearDualVandermonde}, and therefore $\Delta_F \neq 0$.

  For the moreover part, we consider different cases according to $|F|$:
  \begin{itemize}
  \item If $|F| > n$, then $\Delta_{F,g} = 1$.
  \item If $|F| = n$, then $F^\wedge = X^* \cdot x$ is the Chow form of the unique common zero of $F$.
    This is a generic point, so $\Delta_{F,g} = (F,g)^\wedge = \pm g(x) \neq 0$.
  \item If $|F| < n$, then, when calculating $\Delta_{F,g}$ we only consider $g^\sigma$ for $\sigma \subsetneq I_g$.
    Therefore, we cannot distinguish this case from the one where $g$ is also a generic pseudo-split polynomial, so $\Delta_{F,g} \neq 0$ by our main assertion.
  \end{itemize}
  Therefore $\det M_{F,g} = \prod_{G \subseteq F} \Delta_G \Delta_{G,g} \neq 0$.
\end{proof}

\section{A remark regarding the Macaulay resultant}

The previous section yields a curious side result, namely, an expression of the Macaulay resultant of a family of polynomials $(f_i : i \leq n)$ as an inclusion/exclusion expression.
However, this not only requires the polynomials to be indeterminate (this is reasonable, since the resultant is a polynomial in the indeterminate coefficients), but indeed, requires them to be indeterminate pseudo-split polynomials.
In other words, the expression requires the presence of auxiliary indeterminate coefficients for the polynomials $f_i^\sigma$, and specialising those to something ``known'' may cause determinates in both numerator and denominator to vanish.

In this section we propose a different inclusion/exclusion formula for the resultant, which only depends on the indeterminate coefficients of the polynomials.
Like the formula given by Macaulay \cite{Macaulay:Elimination}, it considers products and quotients of determinants of matrices representing, in some sense or another, multiplication by the polynomials $f_i$.
The exclusion/inclusion form of our formula makes is considerably longer than others we found in the literature.
While this makes our formula less useful for actual calculations, it is more symmetric and amenable to formal manipulations.
It was used in an early version of the present paper, but now it is no more than a remark.

Let us fix, as usual, $X = (X_0,\ldots,X_n)$, and let $F = (f_i : i \leq n)$ be an indeterminate family of polynomials, i.e., each $f_i = \sum_\alpha T^*_{i,\alpha} X^\alpha$ is homogeneous of degree $d_i$ with indeterminate coefficients.
As we recalled in \autoref{fct:MacaulayResultant}, the Macaulay resultant $F^\wedge \in \bZ[T^*]$ is irreducible (unless some of the $d_i$ vanish), and uniquely determined if we require that, when specialising $f_i$ to $X_i^{d_i}$, $F^\wedge$ specialise to $1$.

For $\sigma \subseteq n+1 = \{0, \ldots, n\}$, define
\begin{gather*}
  f_\sigma = \prod_{i \in \sigma} f_i,
  \qquad
  X_\sigma = \prod_{i \in \sigma} X_i^{d_i},
  \qquad
  d_\sigma = \deg f_\sigma = \deg X_\sigma = \sum_{i \in \sigma} d_i.
\end{gather*}
Let $\fM_\ell$ denote all monomials of degree $\ell$ in $X$.
For a monomial $\fm$, define
\begin{gather*}
  \rho(\fm) = \bigl\{ i \leq n : X_i^{d_i} \mid \fm \bigr\}.
\end{gather*}
This is the set of $i$ such that we may multiply $\fm$ by $f_i/X_i^{d_i}$ (and keep it a polynomial).
Accordingly, let us define
\begin{gather*}
  [\fm]_F = \frac{\fm f_{\rho(\fm)}}{X_{\rho(\fm)}}.
\end{gather*}
We define $\bA_{F,d}$ to be the $\fM_d \times \fM_d$ matrix, whose $\fm$th row consists of the coefficients of $[\fm]_F$.

\begin{lem}
  \label{lem:MacaulayResultantFactor}
  Assume that $d \geq \sum d_i - n$, and that $d_i \geq 1$ for all $i$.
  Then $F^\wedge \mid \det \bA_{F,d}$ in $\bZ[T^*]$.
\end{lem}
\begin{proof}
  Since $F^\wedge$ is prime, $B = \bZ[T^*]/(F^\wedge)$ is an integral domain, and embeds in an algebraically closed field $L$.
  Let $\varphi\colon \bZ[T^*] \rightarrow L$ be the corresponding morphism.
  Then $\varphi$ sends each $f_i$ to a polynomial $g_i \in L[X]_{d_i}$, and $G^\wedge = \varphi(F^\wedge) = 0$, so the family $G$ admits a common zero $[x] \in \bP^n(L)$.

  Let $V_d(x)$ be the $d$-fold Veronese image of $x$ (i.e., the evaluations of all monomials of degree $d$ at $x$).
  Since $d \geq \sum d_i - n$, each monomial in $\fM_d$ admits a factor $X_i^{d_i}$ for at least one $i$.
  Therefore, $\varphi( \bA_{F,d}) V_d(x) = 0$, so $\varphi\bigl( \det \bA_{F,d} \bigr) = 0$.
  Therefore $\det \bA_{F,d} \in (F^\wedge)$.
\end{proof}

When $n=1$ and $d = d_0 + d_1 - 1$, we have $F^\wedge = \det \bA_{F,d}$ (see \autoref{exm:MacaulayDim1}), but in general there may be other factors to the determinant.
In order to clean them up, so to speak, let us consider variants of the construction.
The general idea consists of adding certain ``exclusions'' to the definition of $[\fm]_F$ and $\bA_{F,d}$.

For $\sigma \subseteq n+1$, define
\begin{gather*}
  [\fm]_F^\sigma = \frac{\fm f_{\rho(\fm) \setminus \sigma}}{X_{\rho(\fm) \setminus \sigma}}.
\end{gather*}
Define $\bA_{F,d}^\sigma$ to be the matrix whose lines are the coefficients of $[\fm]_F^\sigma$ for $\fm \in \fM_{d-d_\sigma}$.
Substituting $X_i^{d_i}$ for each $f_i$ specialises each such matrix to the identity, and its determinant to one.
In particular, the determinant is a non-zero polynomial in $\bZ[T^*]$, and specialises to a non-zero polynomial modulo any prime $p$.
The same holds of other variants we define below.
We shall prove:

\begin{thm}
  \label{thm:MacaulayIncExc}
  Let $d \in \bZ$, and
  \begin{gather*}
    \Xi_{F,d} = \prod_{\sigma \subseteq n+1} \left( \det \bA_{F,d}^\sigma \right)^{(-1)^{|\sigma|}}.
  \end{gather*}
  Then $\Xi_{F,d} \in \bZ[T^*]$.
  If $d \geq \sum_{i \leq n} d_i - n$, then $\Xi_{F,d}$ is the Macaulay resultant $F^\wedge$.
\end{thm}

% \begin{exm}
%   \label{exm:MacaulayDim0}
%   Let $n = 0$, so $F = (f_0)$ and $f_0 = T^*_0 X_0^{d_0}$.
%   The resultant is $T^*_0$.
%   If $d \geq d_0$, then $\bA_{F,d} = T^*_0 I_1$, $\bA_{F,d}^{\{0\}} = I_1$ and $\Xi_{F,d} = T^*_0$.
%   And if $d < d_0$, then $\bA_{F,d}$ is either $I_1$ or $I_0$ (according to whether $d \geq 0$ or not), $\bA_{F,d}^{\{0\}} = I_0$, and $\Xi_{F,d} = 1$.
% \end{exm}

\begin{exm}
  \label{exm:MacaulayDim1}
  Let $n = 1$, so $F = (f_0, f_1)$, and assume that $d = d_0 + d_1 - 1$.
  Then for each monomial $\fm \in \fM_d$ we have either $\rho(\fm) = \{0\}$ or $\rho(\fm) = \{1\}$, but not both.
  The matrix $\bA_{F,d}$ is the Sylvester matrix $S(f_0,f_1)$, $\bA_{F,d}^{\{0\}}$ and $\bA_{F,d}^{\{1\}}$ are the identity matrices in their respective dimensions, and $\bA_{F,d}^{\{0,1\}}$ is empty.

  We obtain the familiar formula for the resultant of two polynomials $f$ and $g$ as $\det S(f, g)$.
\end{exm}

We are going to prove \autoref{thm:MacaulayIncExc} by induction on $n$.
Therefore, assume for the time being that $n \geq 1$.
We are going to treat the index $n$ somewhat differently, so from now on let $\sigma \subseteq n$.

\begin{lem}
  \label{lem:MacaulayExcludeN}
  Let $g_i = f_i(X_0,\ldots,X_{n-1},0) \in \bZ[T^*,X_{<n}]$, and let $G = (g_i : i < n)$.
  Accordingly, for $\sigma \subseteq n$, we have matrices $\bA_{G,d}^\sigma$, indexed by monomials of degree $d$ in $X_{<n}$.
  Then
  \begin{gather*}
    \det \bA_{F,d+d_n}^{\sigma \cup \{n\}}
    = \prod_{d' \leq d} \det \bA_{G,d'}^\sigma.
  \end{gather*}
  For $d' < d_\sigma$ the matrix $ \bA_{G,d'}^\sigma$ is empty, and its determinant is one, so there is no need to specify an explicit lower limit for $d'$.
\end{lem}
\begin{proof}
  Ordering $\fM_d$ as in \autoref{sec:LinearDualVandermonde}, we have
  \begin{gather*}
    \bA_{F,d+d_n}^{\sigma \cup \{n\}} =
    \begin{pmatrix}
      \bA_{G,d}^\sigma & * & * & * \\
                       & \ddots & * & * \\
                       & & \bA_{G,d'}^\sigma & * \\
                       & & & \ddots
    \end{pmatrix}.
  \end{gather*}
  The block $\bA_{G,d'}^\sigma$ correspond to monomials $\fm$ of degree $d-d'$ in $X_n$.
  In particular, if $d < d_\sigma$ then $\bA_{F,d+d_n}^{\sigma \cup \{n\}}$ is empty.
\end{proof}

For an even finer system of exclusions, consider two parameters $\sigma \subseteq n$ and $T \subseteq \cP(n \setminus \sigma)$.
Define
\begin{gather*}
  [\fm]_F^{\sigma,T} =
  \begin{cases}
    [\fm]_F^{\sigma \cup \{n\}} & \rho(\fm) \setminus \sigma \setminus \{n\} \in T
    \\
    [\fm]_F^\sigma & \text{otherwise}.
  \end{cases}
\end{gather*}
Notice that the two cases are distinct only if $n \in \rho(\fm)$.
Define $\bA_{F,d}^{\sigma,T}$ to be the matrix whose lines are the coefficients of $[\fm]_F^{\sigma,T}$ for $\fm \in \fM_{d-d_\sigma}$.
In particular, $\bA_{F,d}^{\sigma,\emptyset} = \bA_{F,d}^\sigma$, and $\bA_{F,d}^{\sigma,\cP(n \setminus \sigma)} = \bA_{F,d + d_n}^{\sigma \cup \{n\}}$.

\begin{lem}
  \label{lem:MacaulayExcludeT}
  Let $\sigma \subseteq n$ and $T' \subseteq T \subseteq \cP(n \setminus \sigma)$.
  Then, in $\bQ(T^*)$, we have
  \begin{gather*}
    \frac{\det \bA_{F,d}^{\sigma,T'}}{\det \bA_{F,d}^{\sigma,T}}
    = \prod_{\tau \in T \setminus T'} \frac{\det \bA_{F,d}^{\sigma \cup \tau}}{\det \bA_{F,d}^{\sigma \cup \tau,\{\emptyset\}}}.
  \end{gather*}
\end{lem}
\begin{proof}
  We may assume that $\sigma = \emptyset$ (replacing $d$ with $d - d_\sigma$, and $f_i$ for $i \in \sigma$ with $1$).
  We may then reduce to the case where $T' = \emptyset$, and from that, to the case where $\tau \in T$ is maximal and $T' = T \setminus \{\tau\}$.
  Then it will suffice to prove that
  \begin{gather*}
    \frac{\det \bA_{F,d}^{\emptyset,T'}}{\det \bA_{F,d}^{\emptyset,T}}
    =
    \frac{\det \bA_{F,d}^\tau}{\det \bA_{F,d}^{\tau,\{\emptyset\}}}.
  \end{gather*}

  Let us order $\fM_d$ so $X_\tau \fM_{d-d_\tau}$ is an initial segment, ordered as $\fM_{d-d_\tau}$.
  Let $\bB$ be the matrix of multiplication of monomials in $\fM_{d-d_\tau}$ by $f_\tau$.
  Let $\fm \in \fM_d$, and let $g$ and $h$ be the polynomials given by the $\fm$th rows of $\bA_{F,d}^T$ and $\bA_{F,d}^{T'}$, respectively:
  \begin{itemize}
  \item If $\fm = X_\tau \fm'$, where $\fm' \in \fM_{d-d_\tau}$, then $g = f_\tau g'$ and $h = f_\tau h'$.
    By maximality of $\tau$ in $T$, it is the only exclusion that possibly applies, so $g'$ is given by the $\fm'$th row of $\bA^{\tau,\{\emptyset\}}_{F,d}$, while $h'$ is given by that of $\bA^\tau_{F,d}$.
  \item Otherwise, $\rho(\fm) \setminus \{n\} \neq \tau$, so $T$ and $T'$ impose the same exclusions, and $g = h$.
  \end{itemize}
  Therefore,
  \begin{gather*}
    \bA_{F,d}^{\emptyset,T}
    =
    \begin{pmatrix}
      \bA_{F,d}^{\tau,\{\emptyset\}} \bB \\
      \bC
    \end{pmatrix},
    \qquad
    \bA_{F,d}^{\emptyset,T'}
    =
    \begin{pmatrix}
      \bA_{F,d}^\tau \bB \\
      \bC
    \end{pmatrix}.
  \end{gather*}
  Let $\bD = \begin{pmatrix} & \bB \\ 0 & & I \end{pmatrix}$.
  As usual, we can specialise it to the identity, so it is invertible over $\bQ(T^*)$.
  Now:
  \begin{gather*}
    \bA_{F,d}^{\emptyset,T} \bD^{-1}
    =
    \begin{pmatrix}
      \bA_{F,d}^{\tau,\{\emptyset\}} & & 0 \\
      & \bC \bD^{-1}
    \end{pmatrix},
    \qquad
    \bA_{F,d}^{\emptyset,T'} \bD^{-1}
    =
    \begin{pmatrix}
      \bA_{F,d}^\tau & & 0 \\
      & \bC \bD^{-1}
    \end{pmatrix}.
  \end{gather*}
  Our assertion follows.
\end{proof}

\begin{lem}
  \label{lem:MacaulayIncExcTwoParts}
  With $n \geq 1$ and $G$ as in \autoref{lem:MacaulayExcludeN}, we have
  \begin{gather*}
    \Xi_{F,d}
    = \frac{\det \bA_{F,d}^{\emptyset, \cP(n) \setminus \{\emptyset\}}}{\det \bA_{F,d+d_n}^{\{n\}}}
    \prod_{k<d_n} \Xi_{G,d-k}.
  \end{gather*}
\end{lem}
\begin{proof}
  By definition,
  \begin{gather*}
    \Xi_{F,d}
    = \prod_{\sigma \subseteq n} \left( \frac{\det \bA_{F,d}^\sigma}{\det \bA_{F,d}^{\sigma \cup \{n\}}} \right)^{(-1)^{|\sigma|}}
    = \prod_{\sigma \subseteq n} \left( \frac{\det \bA_{F,d}^\sigma}{\det \bA_{F,d+d_n}^{\sigma \cup \{n\}}} \right)^{(-1)^{|\sigma|}}
    \prod_{\sigma \subseteq n} \left( \frac{\det \bA_{F,d+d_n}^{\sigma \cup \{n\}}}{\det \bA_{F,d}^{\sigma \cup \{n\}}} \right)^{(-1)^{|\sigma|}}.
  \end{gather*}
  Let us develop the first factor, applying \autoref{lem:MacaulayExcludeT} twice:
  \begin{align*}
    \prod_{\sigma \subseteq n} \left( \frac{\det \bA_{F,d}^\sigma}{\det \bA_{F,d}^{\sigma,\cP(n \setminus \sigma)}} \right)^{(-1)^{|\sigma|}}
    & = \prod_{\sigma \subseteq n} \prod_{\tau \subseteq n \setminus \sigma} \left( \frac{\det \bA_{F,d}^{\sigma\cup\tau}}{\det \bA_{F,d}^{\sigma\cup\tau,\{\emptyset\}}} \right)^{(-1)^{|\sigma|}}
    \\
    & = \frac{\det \bA_{F,d}^\emptyset}{\det \bA_{F,d}^{\emptyset,\{\emptyset\}}}
      = \frac{\det \bA_{F,d}^{\emptyset, \cP(n) \setminus \{\emptyset\}}}{\det \bA_{F,d+d_n}^{\emptyset, \cP(n)}}
      = \frac{\det \bA_{F,d}^{\emptyset, \cP(n) \setminus \{\emptyset\}}}{\det \bA_{F,d+d_n}^{\{n\}}}.
  \end{align*}
  By \autoref{lem:MacaulayExcludeN}, the second factor is
  \begin{align*}
    \prod_{\sigma \subseteq n} \left( \prod_{k<d_n} \det \bA_{G,d-k}^\sigma \right)^{(-1)^{|\sigma|}}
    = \prod_{k<d_n} \Xi_{G,d-k},
  \end{align*}
  completing the proof.
\end{proof}

\begin{proof}[Proof of \autoref{thm:MacaulayIncExc}]
  Let us first show that $\Xi_{F,d} \in \bZ[T^*]$ for all $d \in \bZ$.
  We proceed by induction on $n$.
  If $n = 0$, then the only term in the denominator is $\det \bA_{F,d}^{\{0\}} = 1$.
  For $n \geq 1$, the induction hypothesis together with \autoref{lem:MacaulayIncExcTwoParts} imply that $T^*_n$ does not appear in the denominator of $\Xi_{F,d}$.
  By symmetry, the denominator must be constant.
  Since it divides a polynomial that can be specialised to one, it is one, and $\Xi_{F,d}$ is polynomial.

  By \autoref{lem:MacaulayIncExcTwoParts} we have $\deg_{T^*_n} \Xi_{F,d} = \deg_{T^*_n} \det \bA_{F,d}^{\emptyset, \cP(n) \setminus \{\emptyset\}}$.
  The latter is the number of monomials $\fm = X^\alpha \in \fM_d$ such that $\rho(\fm) = \{n\}$, i.e., such that $\alpha_i < d_i$ for all $i < n$.
  If we add the hypothesis that $d \geq \sum_{i\leq n} d_i - n$, then each such sequence $(\alpha_i : i < n)$ is possible, with $\alpha_n = d - \sum_{i<n} \alpha_i \geq 0$.
  Therefore $\deg_{T^*_n} \Xi_{F,d} = \prod_{i<n} d_i = \deg_{T^*_n} F^\wedge$, and similarly for every other index.

  If $d_i \geq 1$ for all $i$, then $F^\wedge$ is irreducible and divides $\bA_{F,d}$ by \autoref{lem:MacaulayResultantFactor}.
  It cannot divide any other factor of $\Xi_{F,d}$, so $F^\wedge \mid \Xi_{F,d}$ in $\bZ[T^*]$.
  Since the degrees in each $T^*_i$ are the same, $\Xi_{F,d}/F^\wedge$ is a constant in $\bZ$.
  If $d_i = 0$ for some $i$, then $F^\wedge = (T^*_i)^{\prod_{j \neq i} d_j}$, and by the same degree calculation $\Xi_{F,d}/F^\wedge$ is a constant.

  Finally, specialising $f_i$ to $X_i^{d_i}$, we see that $\Xi_{F,d}/F^\wedge = 1$.
\end{proof}

\bibliographystyle{begnac}
\bibliography{begnac}

\end{document}